%% file: brianrig
\documentclass[12pt]{article}
\usepackage[margin=1.4in]{geometry}
\usepackage{amssymb,amsmath}
\usepackage[dvips]{graphics}
\usepackage{epsfig,color}

\usepackage{latexsym}

\newfont{\bb}{msbm10 at 11pt}
\newfont{\bbsmall}{msbm8 at 8pt}

\def\r{\mathbb{R}}

\newcommand{\R}{\mbox{\bb R}}

\newcommand{\C}{\mbox{\bb C}}

\newcommand{\N}{\mbox{\bb N}}

\newcommand{\esf}{\mbox{\bb S}}

\newcommand{\rth}{\R^3}

\def\De{{\Delta}}

\def\ve{{\varepsilon}}

\def\centerbmp#1#2#3{\vskip#2\relax\centerline{\hbox to#1{\special
    {bmp:#3 x=#1, y=#2}\hfil}}}

\newtheorem{theorem}{Theorem}[section]
\newtheorem{lemma}[theorem]{Lemma}
\newtheorem{proposition}[theorem]{Proposition}

\newtheorem{remark}[theorem]{Remark}
\newtheorem{corollary}[theorem]{Corollary}
\newtheorem{definition}[theorem]{Definition}

\newenvironment{proof}{\smallskip\noindent{\it Proof.}\hskip \labelsep}
                          {\hfill\penalty10000\raisebox{-.09em}{$\Box$}\par\medskip}

\begin{document}
\begin{title}
{The number of constant mean curvature isometric immersions of a surface}
\end{title}
\begin{author}
{Brian Smyth \and Giuseppe Tinaglia}
\end{author}
\date{}
\maketitle
\vspace{-1cm}
\begin{center}{\it \small In Memoriam Katsumi Nomizu}\end{center}

\begin{abstract}
In classical surface theory there are but few known examples of surfaces admitting nontrivial isometric deformations and fewer still non-simply-connected ones. 
We consider the isometric deformability question for an immersion $x\colon M \to \rth$ of an oriented non-simply-connected surface with constant mean curvature $H$. We prove that the space of all isometric immersions of $M$ with constant mean curvature $H$ is, modulo congruences of $\rth$, either finite or a circle (Theorem~\ref{main}). When it is a circle then, for the immersion $x$, every cycle in $M$ has vanishing force and, when $H\neq 0$, also vanishing torque. Our work  generalizes a rigidity result for minimal surfaces~\cite{cmw1} to  constant mean curvature surfaces (Theorem~\ref{appl}). Moreover, we identify closed vector-valued 1-forms whose periods give the force and torque. \end{abstract}

\section{Introduction.}

It has long been wondered which smooth surfaces in $\rth$ can undergo a nontrivial isometric deformation and the non-existence of such compact surfaces has been conjectured. Even if we allow immersions, non-compactness, non-completeness or all of these, the occurrence of such deformations for {\it non-simply-connected} surfaces still appears to be rare. Apart from flat surfaces, the only place in classical surface theory where local isometric deformations overtly present themselves is with surfaces of constant mean curvature.  (See also the simply-connected examples of do~Carmo-Dajczer~\cite{doda})

Within the class of isometric immersions of {\it simply-connected} surfaces, every one of constant mean curvature admits a canonical $2\pi$-periodic isometric deformation --- the associate deformation (see \S~\ref{deformation}) --- through immersions with the same mean curvature; this family also captures (to within congruences of $\rth$, i.e. orientation-preserving isometries of $\rth$) all isometric immersions of that surface with this constant mean curvature. However once the constant mean curvature surface $M$ has topology (i.e. $\pi_1(M)\neq 0$) there is of necessity a continuum of latent period conditions which must be satisfied for the canonical deformation, at the level of the universal cover $\widetilde{M}$,  to descend to $M$. Our purpose here is to give necessary conditions for the associate isometric deformation to exist.

A countable set of invariants associated with a constant mean curvature immersion $x$ arises from two naturally defined closed vector-valued 1-forms $\omega$ and $\sigma$ on $M$, introduced in \S~\ref{force}, called here the {\it force} and the {\it torque} forms. Thus their periods over a cycle $\gamma$ in $M$, $W([\gamma])=\int_{\gamma}\omega$ and $T([\gamma])=\int_{\gamma}\sigma$, depend only on the homology class of $\gamma$ and are called the force and torque of that class for the immersion $x$  (see also \cite{kks1,ku2,loraf1}).

\begin{theorem}\label{main}

Let $x\colon M\to \rth$ be an isometric immersion of a smooth oriented surface with constant mean curvature $H$. To within congruences, the family of all isometric immersions of $M$ with constant mean curvature $H$ is either finite or a circle.

If the family is a circle then for the immersion $x$ every cycle has vanishing force and, when $H \neq 0$, vanishing torque also.
\end{theorem}

There was already a class of results on the isometric indeformability of minimal or constant mean curvature surfaces with topology (see for instance \cite{cmw1, ku2, me6, mr3, mt2, perez4}). Typically, these follow for us from Theorem \ref{main} by exhibiting a cycle with nonzero force. In fact, in \S~5, we discuss geometric conditions which guarantee the existence of such a cycle. The central questions remaining are:

\begin{quote}
\it Are there non-simply-connected constant mean curvature surfaces satisfying the vanishing force and torque conditions of Theorem~\ref{main} and do these conditions then guarantee the existence of a constant mean curvature isometric deformation? 
\end{quote}

Both questions are known to be answered affirmatively for minimal surfaces (see~\S~\ref{force}). However, for non-simply-connected surfaces we know of no constant mean curvature isometric deformation when $H\neq 0$. Franz Pedit informs us that there are immersed constant mean curvature cylinders with vanishing force for which the associate deformation does not exist; the torque is unknown. Perhaps the best problem is: 

\begin{quote}
\it A compact constant mean curvature surface admits no non-trivial constant mean curvature isometric deformations.
\end{quote}

Among complete Riemannian surfaces $(M,g)$ the flat ones $\r^2$ and $\esf^1\times \r$ can be isometrically immersed in $\rth$ with any constant mean curvature $H$ (as circular cylinders of radius $\frac{1}{2|H|}$). We note that if the metric is not flat and $x$ is an isometric immersion of $(M,g)$ with constant mean curvature $H$ then $H^2$ is uniquely determined by the metric $g$: $-x$ is an isometric immersion with constant mean curvature~$-H$.

\noindent {\bf Acknowledgements.} We are grateful to Rob Kusner and Bill Meeks for their interest and helpful remarks on the results of this paper. 

\section{Isometric deformation of surfaces}

Let $x\colon M \to \rth$ be an immersion of a smooth oriented surface. The differential $x_*$ of $x$ is given by $x_*(X)=Xx$ where the right-hand side is the derivative of the vector-valued function $x$ with respect to $X$. The induced metric $g$ is given by
$$g(X,Y)=\langle x_* (X),x_*(Y) \rangle = \langle Xx, Yx \rangle $$
 where $\langle\, , \rangle$ is the Euclidean metric on $\rth$. Let $J$ denote the complex structure induced by  the orientation of $M$ and let $\xi$ be the oriented unit normal field to the immersion $x$. The second fundamental form $A$ of $x$ is defined by 
$$X\xi=-x_*(AX).$$
The Gauss equation is 
$$\text{det}A=K$$ 
where $K$ is the curvature of the metric $g$ and Codazzi's equation is 
$$(\nabla_XA)Y=(\nabla_YA)X$$
where $\nabla$ is the Levi-Civita connection of the metric $g$; from now on the metric will be denoted by $\langle\, ,\rangle$. These equations come from differentiating the structure equation
\begin{equation}\label{structure}XYx=x_*(\nabla_XY)+\langle AX,Y \rangle \xi.\end{equation}
Now suppose $x_t\colon M \to \rth$ is a smooth 1-parameter family of immersions each inducing the same metric $\langle\, ,\rangle$; this is called an {\it isometric deformation}. The unit normal field and second fundamental form of each immersion $x_t$ are denoted by $\xi_t$ and $A_t$ respectively. From now on prime denotes differentiation with respect to $t$.

Since $\langle Xx_t,Xx_t \rangle$ is independent of $t$ we know that $\langle Xx'_t,Xx_t\rangle=0$. Hence $Xx'_t=k(x_t)_*(JX)+p(X)\xi_t$ where $k$ is a function on $M$ and $p$ is a 1-form on $M$, both dependent on $t$.

Since $\langle \xi_t,\xi_t \rangle=1$ it follows that $\xi'_t=(x_t)_*(Z)$ where $Z$ is a vector field on $M$, dependent on $t$. Since $\langle Xx_t,\xi_t\rangle=0$ it follows 
$$\langle Xx'_t,\xi_t \rangle+\langle Xx_t , \xi'_t \rangle=0$$
and so $p(X)=-\langle X, Z \rangle$. Thus 
\begin{equation}\label{derivative}
Xx'_t=k\text{ }(x_t)_*(JX)-\langle X, Z \rangle \xi_t.
\end{equation}
In continuing the computation for the deformation we will drop the subscript $t$. From equation \eqref{structure}
\begin{align}
XYx'= & X(kx_*(JY)-\langle Y, Z \rangle \xi) \notag \\ 
= & X(k)x_*(JY)+k[x_*(\nabla_X JY) +\langle AX,JY \rangle \xi]- X\langle Y,Z \rangle\xi+\langle Y,Z \rangle x_*(AX) \notag\\ 
= & x_*(X(k)JY+kJ\nabla_XY +\langle Y, Z\rangle AX)+(k\langle AX, JY\rangle - X\langle Y, Z \rangle )\xi. \label{secondderiv} 
\end{align}

Differentiating equation \eqref{structure} with respect to $t$ and using equation \eqref{derivative} gives
\begin{align}
XYx'=& \nabla_XY(x')+\langle A'X,Y\rangle \xi +\langle AX, Y \rangle x_*(Z) \notag \\
=& kx_*(J\nabla_XY)-\langle \nabla_XY,Z\rangle\xi +\langle A'X,Y \rangle \xi +\langle AX,Y\rangle x_*(Z) \notag \\
=& x_*(kJ\nabla_XY+\langle AX,Y\rangle Z)+(\langle A'X,Y\rangle -\langle \nabla_XY,Z\rangle)\xi. \label{four}
\end{align}

Comparing normal components in equations \eqref{secondderiv} and \eqref{four} we obtain
\begin{equation}\label{one} A'X=-kJAX-\nabla_XZ \end{equation}
and comparing tangential components
$$X(k)JY=\langle AX,Y\rangle Z- \langle Y,Z\rangle AX $$
which is equivalent to
\begin{equation}\label{two}\nabla k=-AJZ.\end{equation}

Equations \eqref{one} and \eqref{two} are the integrability conditions of the deformation. Finally note
\begin{equation}\label{three}
Xx'=(x_*(JZ)+k\xi)\times Xx_t
\end{equation}
where $\times$ denotes the cross product on $\rth$. The vector-valued function $\eta=x_*(JZ)+k\xi$ is called the Drehriss~\cite{blas} and is used by Burns and Clancy in their recent work~\cite{bcl}.

\section{The associate surfaces of a constant mean curvature surface}\label{deformation}

Let $x \colon M \to \rth$ be an oriented surface with constant mean curvature $H$ then, for each $t\in [0,2\pi]$, the symmetric tensor field
\begin{equation}\label{eight}
A_t=\cos (t) (A-HI) +\sin (t) J(A-HI)+HI
\end{equation}
 satisfies the Gauss-Codazzi equations with respect to the induced metric $\langle\, , \rangle$ and $\frac{\text{Tr}A_t}{2}\equiv H$. If $M$ is simply-connected then, by the fundamental theorem of surface theory, we obtain a one-parameter family of isometric immersions $x_t$  with second fundamental form $A_t$ and therefore constant mean curvature $H$ and $x_0=x$. These immersions are uniquely determined to within a congruence or rigid motion of $\mathbb{R}^3$, that is, an orientation preserving isometry of $\rth$. Without loss of generality, we may assume $x(p_0)=0$ for some $p_0\in M$ and normalize the family $x_t$ by requiring that, for all $t$, $x_t(p_0)=0$, $\xi_t(p_0)=\xi(p_0)$ and $(x_t)_{*_{p_0}}=x_{*_{p_0}}\circ \text{Rot}_{p_0}(-t)$ where $\text{Rot}_{p_0}(\theta)$ denotes the oriented rotation of the tangent plane $T_{p_0}M$ (with the induced metric) through an angle $\theta$.
The resulting normalized isometric deformation $x_t \colon M \to \rth$, $t\in[0,2\pi]$ is called the  {\it associate deformation} here (see also~\cite{Bonnet, ca1, la3}). 

For example, let $x \colon M \to \rth$ be an oriented {\it simply-connected minimal} surface and choose the origin of $\rth$ on the surface, i.e. $x(p_0)=0$ for a certain $p_0\in M$. Since $\De x \equiv 0$, where $\De$ is the Laplace operator of the induced metric on $M$, we have a complex conjugate $y\colon M\to \rth$ of $x$, unique to within a translation of $\rth$. We may therefore assume $y(p_0)=0$ also. Then
$$x_t=\cos(t)x+\sin(t)y $$
is a 1-parameter family of minimal isometric immersions of $M$ into $\rth$ with second fundamental form $A_t$ as given in equation~\eqref{eight} above. Since by the Cauchy-Riemann equations $x_*(X)=y_*(JX)$ and $x_*(JX)=-y_*(X)$ it is easy to see that $\xi_t$ does not change with $t$ and $(x_t)_{*_p}=x_{*_p}\circ \text{Rot}_p(-t))$.

Returning now to the constant mean curvature case, if $M$ is not simply-connected then we may lift $x$ to $\widetilde{x}\colon \widetilde{M}\to\rth$ and let $\widetilde{A}_t$ denote the lift of $A_t$ to $\widetilde{M}$. By the earlier discussion, $\widetilde{A}_t$ is the second fundamental form of an isometric immersion $\widetilde{x}_t \colon \widetilde{M}\to \rth$ with constant mean curvature $H$ and is unique to within a motion of $\rth$. Fixing $p_0\in M$ and $\widetilde{p}_0\in\widetilde{M}$ over $p_0$, we may assume $\widetilde{x}_t(\widetilde{p}_0)=0$, $\widetilde{\xi}_t (\widetilde{p}_0)=\xi(p_0)$ and $(\widetilde{x}_t)_{*_{\widetilde{p}_0}}=\widetilde{x}_{*_{\widetilde{p}_0}}\circ \text{Rot}_{\widetilde{p}_0}(-t)$ for all $t$. With this normalization we obtain a smooth isometric deformation  
$$\widetilde{x}_t \colon \widetilde{M}\to \rth, \quad 0\leq t\leq2\pi $$ with constant mean curvature $H$.
Of course $\widetilde{x}_0$ projects to $x \colon M \to \rth$. Let $S=\{t\in [0,2\pi) \mid \widetilde{x}_t \text{ projects to } x_t\colon M \to \rth \};$ this set of immersions of $M$ in $\rth$ will be called the {\it associate family} for $x\colon M \to \rth$. 

\begin{lemma}\label{associate}
If $x_m \colon M \to \rth$, $m=1,2$, are isometric immersions with constant mean curvature $H$ then $x_2$ is congruent to a unique associate of $x_1$.

Thus if $x \colon M \to \rth$ is an immersion with constant mean curvature $H$, all other isometric immersions of $M$ in $\rth$ with the same constant mean curvature  occur, to within congruences, in the family of associates of $x$. 
\end{lemma}

\begin{proof}
Locally on $M$ we may choose positive isothermal coordinates $(u,v)$, i.e. $\{\frac{\partial}{\partial u},\frac{\partial}{\partial v}\}$ is a positively oriented frame and the metric is of the form $\langle \, ,\rangle=e^{2\rho}(du^2+dv^2)$. The second fundamental form of $x_m$ with respect to the coordinate frame is written
$$A_m=\left[
\begin{array}{cc}
H+\alpha_m & \beta_m\\
\beta_m &H-\alpha_m
\end{array} 
\right]$$
Now, $\omega_m=\langle A_m \frac{\partial}{\partial w},\frac{\partial}{\partial w}\rangle$, is a complex function in the coordinate $w=u+iv$, where $\frac{\partial}{\partial w}=\frac12(\frac{\partial}{\partial u}-i\frac{\partial}{\partial v})$, and $\langle \, , \rangle$ and $A$ are extended by complex linearity. Clearly, $\Omega_m=\omega_m dw^2$ is a well defined complex quadratic differential on $M$; a simple computation gives $\omega_m=-\frac{i}{2}e^{2\rho}(\beta_m+i\alpha_m)$. Codazzi's equation in these isothermal coordinates reads $$(\omega_m)_{\overline{w}}=\frac12e^{2\rho}H_w$$ and, since $H$ is constant, $\omega_m$ is holomorphic. 

%

Since $|\omega_m|^2=\frac{e^{4\rho}}{4}(\beta_m^2+\alpha_m^2)=\frac{e^{4\rho}}{4}(H^2-K)$, by Gauss' equation, unless $x \colon M \to \rth$ is a round sphere the umbilical points are isolated and it follows that the meromorphic function $\frac{\omega_2}{\omega_1}$ is constant and of modulus one. Hence, $\omega_2=e^{-it}\omega_1$ for some real number $t$. It follows easily that $$A_2=\cos (t) (A_1-HI) +\sin (t) J(A_1-HI)+HI.$$
Hence, by the fundamental theorem of surface theory $x_2$ is congruent to an associate of $x_1$.
\end{proof}


We now consider the structure of $S$ for a constant mean curvature surface $M$ with topology. The first and most interesting question is whether $S$ contains an open interval.
Assume $S$ contains an interval $[0,\ve)$ then we have the associate deformation
$$x_t \colon M \to \rth $$ 
for $0\leq t\leq \ve$. Since, by \eqref{eight},  $A'=J(A-HI)$ the integrability condition \eqref{one} becomes 
\begin{equation}\label{seven}
\nabla_XJZ=(k+1)AX-HX.
\end{equation}
Replacing $Y$ by $JZ$ in the structure equation \eqref{structure} and using equation~\eqref{two} and~\eqref{seven}
\begin{align*}
Xx_*(JZ)&= x_*(\nabla_XJZ)+\langle AX,JZ\rangle \xi\\
&=x_*((k+1)AX-HX)+\langle X,AJZ\rangle \xi\\
&=x_*((k+1)AX-HX)-\langle X,\nabla k\rangle \xi\\
&=-X((k+1)\xi+Hx)
\end{align*}
so that $(x_t)_*(JZ)+(k+1)\xi_t+Hx_t=V_t$ is a constant vector field along each surface in the variation. From the normalization in the definition of the associate deformation we obtain $Z(p_0)=0$ and $(x_t)_{*_{p_0}}X=x_{*_{p_0}}\circ \text{Rot}_{p_0}(-t)X$ gives 
$$\frac{d}{dt}(x_t)_{*_{p_0}}X=x_{*_{p_0}}(-\sin(t)X-\cos(t)JX)=-x_{*_{p_0}}\circ \text{Rot}_{p_0}(-t)JX.$$
On the other hand, from equation \eqref{derivative},
$$\frac{d}{dt}(x_t)_{*_{p_0}}X=X(x'_t)\mid_{p_0}=k(x_t)_{*_{p_0}}(JX)
-\langle X,Z\rangle \xi_t(p_0)=k(p_0)x_{*_{p_0}}\circ \text{Rot}_{p_0}(-t)JX.$$
Hence $k(p_0)=-1$ and so $V_t(p_0)=0$ and therefore $V_t\equiv 0$, i.e. $$(x_t)_*(JZ)+(k+1)\xi_t+Hx_t\equiv 0.$$
Now 
$$Xx'_t=((x_t)_*(JZ)+k\xi_t)\times Xx_t= -(Hx_t+\xi_t)\times Xx_t$$
and since $\langle Xx_t, \xi_t \rangle=0$ it follows, also on differentiating with respect to $t$, that 
$$\xi'_t=-(Hx_t+\xi_t)\times \xi_t.$$
We collect these facts in the following lemma
\begin{lemma}
If the associate deformation of $x\colon M \to \rth$ exists then
$$Xx'=-(Hx+\xi)\times Xx,$$
$$\xi'=-(Hx+\xi)\times \xi.$$
\end{lemma}

\begin{lemma}\label{newlemma}
Let $M$ be a complete oriented Riemannian 2-manifold and $x_1$ and $x_2$  isometric immersions with constant mean curvature $H_1$ and $H_2$. Either
\begin{itemize}
\item[(i)] $M$ is flat and each $x_i$ is a circular cylinder of radius $\frac{1}{2|H_i|}$, or
\item[(ii)] $M$ is not flat, $H_1^2=H_2^2$ and $x_2$ is an associate of $\pm x_1$.
\end{itemize}
\end{lemma}
\begin{proof}
If $H_2=\pm H_1$ then {\it (ii)} follows from Lemma~\ref{associate}. If $H_1^2<H_2^2$ then, since $H_1^2-K\geq 0$ on $M$, $H_2^2-K$ is positive and bounded away from zero on $M$. The complete metric $g_0=\sqrt{H_2^2-K}g$ on $M$ is flat. This is because if $g=e^{2\rho}(du^2+dv^2)$ and $\Omega_2=\omega_2dw^2$ are local representation of the metric and the Hopf differential then $g_0=|\omega_2|(du^2+dv^2)$ is the local representations of $g_0$; this is flat since $\omega_2$ is holomorphic (see \S~\ref{deformation}). The universal cover of $(M,g_0)$ is conformally $\mathbb{C}$ and $\frac{\Omega_1}{\Omega_2}$ is a holomorphic function on $\mathbb{C}$ with $\frac{|\Omega_1|}{|\Omega_2|}=\frac{\sqrt{H_1^2-K}}{\sqrt{H_2^2-K}}<1$ and therefore constant. Thus $K$ is constant and for conformal reasons it must be zero. It follows easily that $x_i$ is an immersion of $\mathbb{R}^2$, or a cylindrical quotient thereof, as a cylinder of radius $\frac{1}{2|H_i|}$ in $\rth$. 

\end{proof}
\section{The force and torque 1-forms on an immersed surface of constant mean curvature}\label{force}

To motivate the notions of force and torque take an embedded oriented surface 
$$x\colon M \to \rth $$
with oriented normal $\xi$ and constant mean curvature $H\geq 0$. Imagining the surface as a liquid membrane in equilibrium under a constant normal pressure field $F$, the equilibrium equation is $F=-2H\tau \xi$ \cite{smy2}, where $\tau$ is the surface tension of the membrane. We may assume $\tau=1$.

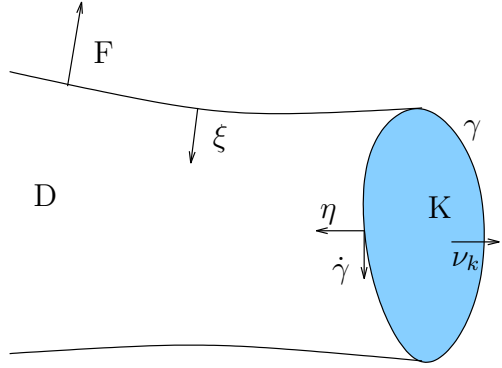
\begin{figure}[htbp]
\begin{center}    
\begin{minipage}[t]{0.5\textwidth}
   \input{orientation.pstex_t}
  \caption{Force and Torque}\label{orientation}
    \end{minipage}
\end{center}
\end{figure}

Take a compact domain $D$ in $M$ and along each oriented boundary component $\gamma$ we insert a smooth embedded cap  
$$k \colon K \to \rth $$
that is $k(\partial K)=x(\gamma)$. The orientation of $M$ determines the orientation of each connected component of $\gamma$: in fact, let $\eta=J\dot{\gamma}$ be the oriented unit normal to $\gamma$ in $M$ then $\xi=x_*\dot{\gamma}\times x_*\eta$. In turns, the orientation of $\gamma$ determines an orientation on $K$ and we let $\nu_k$ be the oriented unit normal to $K$ (see Figure~\ref{orientation}). 

Considering the domain $D$ with caps inserted on each boundary component $\gamma$, the resulting closed surface is maintained in equilibrium by the application of a total restorative force on each end --- to counter the inherent forces due to pressure and surface tension of that end --- which total
$$2H\int_K\nu_k da_k+\int_\gamma \eta ds $$
where $da_k$ also denotes the area element in $K$.

Define $W([\gamma])=2H\int_K\nu_k da_k+\int_\gamma \eta ds $ as the force of the component $\gamma$.

Let $\omega_0$ be the 1-form on $M$ defined by $\omega_0(X)=Hx\times x_*(X)$ then $d\omega_0=2H\xi da$ where $da$ is the area element of $M$. The corresponding 1-form $\omega_0^k$ on $K$, defined by $\omega_0^k(X)=Hk\times k_*(X)$, satisfies $d\omega_0^k=2H\nu_k da_k$. If $\omega$ is the 1-form defined on $M$ by $\omega(X)=(Hx+\xi)\times x_*(X)$ then
$$\int_\gamma \omega= \int_\gamma \omega_0+\int_\gamma\eta ds = \int_\gamma\omega_0^k+\int_\gamma\eta ds,$$
since $\omega_0^k=\omega_0$ along $\gamma$. By Stokes' theorem
\begin{equation}\label{cap} 
\int_\gamma\omega=\int_K d\omega_0^k+\int_\gamma\eta ds= 2H\int_K \nu_k da_k+\int_\gamma\eta ds=W(\gamma).
\end{equation}

Now $\omega$ is easily checked to be a closed 1-form on $M$ for any immersed oriented surface
$$x\colon M\to\rth$$
of constant mean curvature $H$; we call it the {\it force form}. Thus the quantity
$$W([\gamma])=\int_\gamma \omega$$
depends only on the homology class of the cycle $\gamma$. This will be called the {\it force} of the cycle $\gamma$ for the immersion $x$.

Returning again to the domain $D$ with ends capped as above the torque of the inherent forces of pressure and surface tension at the end $\gamma$ totals
$$2H\int_K k\times \nu_k da_k+\int_\gamma x\times \eta ds.$$

Define the torque of $\gamma$ by $$ T(\gamma)=2H\int_K k \times \nu_k da_k+\int_\gamma x \times \eta ds.$$ 
 
Define $\sigma_0(X)=\frac23 Hx \times ( x\times x_*(X)).$ We can easily compute $d\sigma_0=2H x\times \xi da$. The corresponding 1-form $\sigma_0^k$ on $K$ defined by $\sigma_0^k(X)=\frac23 Hk \times (k \times k_*(X))$ satisfies $d\sigma_0^k=2Hk \times \nu_k da_k$. Thus
$$2H\int_K k\times \nu_k da_k=\int_{\partial K}\sigma_0^k=\int_\gamma \sigma_0$$
since $\sigma_0^k=\sigma_0$ along $\gamma$. Then $T(\gamma)=\int_\gamma\sigma_0+\int_\gamma x \times \eta ds =\int_\gamma \sigma$ where $\sigma$ is the $1$-form defined by
$$\sigma(X)=\frac23 H x\times (x \times x_*(X)) + x \times x_*(JX)= \frac{1}{3}x\times [2(Hx+\xi)\times x_*(X)+x_*(JX)].$$
Again it is easy to check that $\sigma$ is a closed 1-form on any immersed surface of constant mean curvature and we call it the {\it torque form}. Now 
$$T([\gamma])=\int_\gamma \sigma$$
depends only on the homology class of $\gamma$ and is called the {\it torque} of the cycle $\gamma$ for the immersion $x$.

To prove the second part of Theorem~\ref{main} announced in the introduction we need only show

\begin{lemma}\label{mainlemma}
Let $x \colon M\to \rth$ be an immersed surface of constant mean curvature $H$ admitting a nontrivial isometric deformation through surfaces of constant mean curvature $H$ then for the immersion $x$
\begin{enumerate}
\item[(i)] $\omega$ is exact;
\item[(ii)] $\sigma$ is exact if $H\neq 0$.
\end{enumerate}
\end{lemma}

\begin{proof}
Consider the associate family $x_t\colon M\to \rth$ with $x_0=x$ then, from the assumption of the lemma, $x_t$ is defined for $0\leq t <\ve$, for some $\ve>0$. As we saw in \S 3,
$$Xx'=-(Hx+\xi)\times Xx = -\omega (X).$$ 
Hence $\omega$ is exact.

We begin by calculating $Xx''$. Write $P=Hx+\xi$
\begin{align*}
Xx''&=-(P\times Xx)'=-P'\times Xx + P\times(P\times Xx)\\
&=-Hx'\times Xx-\xi'\times Xx+P\times(P\times Xx)\\ 
&= -HX(x'\times x)+HXx'\times x-\xi'\times Xx+P\times(P\times Xx)\\
&= -HX(x'\times x)-H(P\times Xx)\times x +(P\times \xi)\times Xx +P\times (P\times Xx)
\end{align*}
\begin{align*}
X(x''+&Hx'\times x)= (Hx+P)\times (P\times Xx)+(P\times \xi)\times Xx\\
&=(2Hx+\xi)\times(P\times Xx)-Xx\times(P\times \xi)\\
&= 2Hx\times (P\times Xx)+ \xi\times(P\times Xx)-Xx\times(P\times\xi)\\
&=2Hx\times(P\times Xx)-P\times(Xx\times \xi)-Xx\times(\xi\times P)-Xx\times(P\times \xi)\\
&=2Hx\times(P\times Xx)+(Hx+\xi)\times x_*(JX)
\end{align*}
Thus
$$X(x''+Hx'\times x +x)=Hx\times\{2P\times Xx+x_*(JX)\}=3H\sigma(X).$$
Hence, when $H\neq 0$, $\sigma$ is exact if there exists an isometric deformation of $x$. \end{proof}

To complete the proof of Theorem~\ref{main} it remains to show:

\begin{lemma}\label{last}
Let $x \colon M\to \rth$ be an immersion of a smooth oriented surface with constant mean curvature $H$. Then either

\begin{itemize}

\item[(i)] to within congruences there are only finitely many isometric immersions of constant mean curvature $H$, or

\item[(ii)] the associates $x_t \colon M \to \rth$ exist for all $t\in[0,2\pi]$.

\end{itemize}
\end{lemma}

\begin{proof}
Let $\widetilde{M}$ be the universal cover of $M$ with the lifted metric and complex structure (denoted $\langle\, , \rangle$ and $J$ respectively), $\pi \colon \widetilde{M} \to M$ the projection and $\mathcal{D}$ be the group of deck transformations of this cover which are, of course, orientation preserving isometries. If $\widetilde{A}$ is the lift to $\widetilde{M}$ of the second fundamental form $A$ of the immersion $x$ then $\sigma_{*_p}\widetilde{A}(p)(\sigma_{*_p})^{-1}=\widetilde{A}(\sigma(p))$ for all $\sigma \in\mathcal{D}$ and all $p\in \widetilde{M}$. The lift $\widetilde{A}_t$ of 
$$A_t=\cos t(A-HI)+\sin t J(A-HI)+HI$$
to the universal cover is the second fundamental form of the associate family
$$\widetilde{x}_t \colon \widetilde{M}\to \rth $$
defined in \S~\ref{deformation}. Since each deck transformation $\sigma\in \mathcal{D}$ commutes with $\widetilde{A}$ and $J$, it preserves the lifted second fundamental form $\widetilde{A}_t$, that is, 
$$\sigma_{*_p}\widetilde{A}_t(p)(\sigma_{*_p})^{-1}=\widetilde{A}_t(\sigma(p))$$
for all $p\in\widetilde{M}$. It follows that $\widetilde{x}_t\circ\sigma$ and $\widetilde{x}_t$ have the same second fundamental form
$$(\sigma_{*_p})^{-1}\widetilde{A}_t(\sigma(p))\sigma_{*_p}=\widetilde{A}_t(p)$$
at each $p$. Hence $\widetilde{x}_t\circ\sigma=\Phi_t(\sigma)\circ \widetilde{x}_t$, where $\Phi_t(\sigma)\in \mathcal{M}$ the group of motions of $\rth$. It easy to see that
$$\Phi_t \colon \mathcal{D} \to \mathcal{M} $$
is a homomorphism for each $t\in[0,2\pi]$ and $\widetilde{x}_t$ projects to $x_t \colon M\to \rth$ if and only if $\Phi_t(\mathcal{D})=\{I\}$.

Let $S=\{t\in [0,2\pi]\mid \widetilde{x}_t \text{ projects to } M\}$. Assuming $S$ is infinite there exists an infinite sequence of points $\{t_n\}$ of points in $S$ which we may assume converges to some $t_0\in[0,2\pi]$. Since $\Phi_{t_n}(\sigma)=I$ for all $n$ it follows, by continuity, that $\Phi_{t_0}(\sigma)=I$ for all $\sigma\in D$. So $t_0\in S$. In the corresponding matrices $(\Phi_{t_n}(\sigma))_{ij}=\delta_{ij}$ for each $t_n$ and so, by the Mean Value Theorem, $(\Phi'_{s_n}(\sigma))_{ij}=0$ for a sequence $s_n$ converging to $t_0$. Hence $\Phi'_{t_0}(\sigma)=0$ for all $\sigma\in D$. Repeating the argument gives that the derivatives $\Phi^{(n)}_{t_0}(\sigma)=0$ for all $\sigma\in D$ and all integers $n\geq 1$. Since, for each $\sigma \in \mathcal{D}$, $\Phi_t(\sigma)$ is an analytic curve in $\mathcal{M}$, it follows that $\Phi_t(\sigma)\equiv I$ for all $t$ and all $\sigma \in \mathcal{D}$. Hence $S=[0,2\pi]$ if $S$ is infinite. This proves Lemma~\ref{last}.
\end{proof}

\begin{remark}\label{minimal}
For minimal immersions the associate $x_{\pi}$ is in fact $-x$ and thus always exists. However, if there exists an associate $x_s$, $s\notin\{0,\pi\}$, then they all exist and $$x_t=\frac{1}{\sin s}\{\sin (s-t)x + \sin (t) x_s\}.$$ 
\end{remark}

Regarding the finiteness part of Theorem~\ref{main}, we note Meeks' conjecture that for a complete isometric embedding $x \colon M \to \rth$ with constant mean curvature $H$ there are no other isometric immersions with the same mean curvature $H (\neq 0)$: when $H=0$ and $x$ is not a helicoid any other minimal isometric immersion is congruent to $-x$ (see~\cite{me14,mpe2}).

As we see in the next result the existence of associates of a constant mean curvature immersion $x \colon M \to \rth$ is guaranteed if there are isometries of the induced metric on $M$ which do not extend under the immersion $x$ to a congruence of $\rth$. However, we know of no non-simply-connected examples with nonzero constant mean curvature where an associate exists except for certain nodoids $x$, i.e. non-embedded Delaunay surfaces, and then only for $x_{\pi}$; for unduloids, i.e. embedded Delaunay surfaces, no associate exists (see for instance \cite{mt2}).

\begin{proposition}
Let $x\colon M \to \rth$ be an immersion of a smooth oriented surface with constant mean curvature $H$.
Let $I(M)$ be the group of orientation preserving isometries of $M$ and let $I_0(M)$ be the subgroup of such isometries extending under $x$ to a congruence of $\rth$. Then $I_0(M)$ extends under each associate $x_t$ of $x$ to a group of isometries of $\rth$. If $I(M)\backslash I_0(M)$ is infinite then $x$ is isometrically deformable through immersions of constant mean curvature. If $I(M)\backslash I_0(M)$ is finite of order $m$ then each of the associates $x_{\frac{2\pi r}{m}}$ exists, $r=1,...,m$.



\end{proposition}

\begin{proof}
If $\phi\in I_0(M)$ then $\phi_{*_p}A_p (\phi_{*_p})^{-1}=A_{\phi(p)}$ for all $p\in M$. Since $\phi$ preserves the metric and complex structure of $M$ we have $\phi_{*_p}A^t_p (\phi_{*_p})^{-1}=A^t_{\phi(p)}$ for all $p\in M$. Hence $\phi$ extends under $x_t$ to an isometry $\Phi$ of $\rth$.

Let $\sigma\in I(M)$. Comparing the isometric immersions $x$ and $y=x\circ \sigma$ the respective oriented normals to these maps at $p$ are $\xi(p)$ and $N(p)=\xi(\sigma(p))$. If $B$ is the second fundamental form of $y$ then, by its definition, $(XN)_p=-y_{*_p}(B(p)X)=-x_{*_{\sigma(p)}}\sigma_{*_p}(B(p)X)$. But, since $N=\xi\circ \sigma$, $$(XN)_p=\xi_{*_{\sigma(p)}}(\sigma_{*_p}X)=-x_{*_{\sigma(p)}}(A(\sigma(p))\sigma_{*_p}X).$$ Taken together, these give
$$B(p)=(\sigma_{*_p})^{-1}A(\sigma(p))\sigma_{*_p},$$
so $y$ has constant mean curvature $H$. By Lemma \ref{associate}, $y=x\circ\sigma$ is congruent to a unique associate $x_{t(\sigma)}$ of $x$. This defines a map $t\colon I(M)\to [0,2\pi)$. Obviously, $t(\sigma)=0$ if and only if $\sigma\in I_0(M)$. To complete the proof of the theorem we must show that $t$ is a homomorphism.

For $\sigma,\tau\in I(M)$ let $C$ denote the second fundamental form of $x\circ\sigma\circ\tau$. As before, 
\begin{align*}
C(p)&=((\sigma\circ\tau)_{*_p})^{-1}A((\sigma\circ\tau)(p))(\sigma\circ\tau)_{*_p}\\
&=(\tau_{*_p})^{-1}(\sigma_{*_{\tau(p)}})^{-1}A(\sigma(\tau(p)))\sigma_{*_{\tau(p)}}\tau_{*_p}\\
&=(\tau_{*_p})^{-1}B(\tau(p))\tau_{*_p}
\end{align*}

Since $x\circ \sigma$ is congruent to the associate $x_{t(\sigma)}$ of $x$ we have $B=e^{it(\sigma)}(A-HI)+HI$, where $e^{it}(A-HI)=\cos t(A-HI) +\sin t J(A-HI)$. Thus
\begin{align*}
C(p)&=e^{it(\sigma)}(\tau_{*_p})^{-1}(A-HI)(\tau(p))\tau_{*_p}+HI(p)\\
&=e^{it(\sigma)}[(\tau_{*_p})^{-1}A(\tau(p))\tau_{*_p}-HI(p))]+HI(p)\\
&=e^{it(\sigma)}\{e^{it(\tau)}(A-HI)(p)\}+HI(p),
\end{align*}
since $x\circ\tau$ is congruent to $x_{t(\tau)}$. Thus 
$$C(p)=e^{i(t(\sigma)+t(\tau))}(A-HI(p))+HI(p).$$
and hence $t(\sigma\circ\tau)=(t(\sigma)+t(\tau)) \text{ mod } 2\pi$ and the map $t$ defines a homomorphism.
\end{proof}

\section{Applications}

We begin with the observation that there exist complete, immersed minimal surfaces with genus zero and finitely many ends admitting a nontrivial isometric deformation through minimal surfaces. First, let us recall the Weierstrass representation theorem, see for instance \cite{os1},                                                                   

\begin{proposition}\label{wt}
Let $x \colon M \to \rth$ be a conformal minimal immersion of a Riemann surface $M$. Let $g$ be the
stereographic projection of its Gauss map, $dh=dx_3-idx_3\circ J$; $g$ and $dh$ are holomorphic. Then (up to a translation) 

\begin{equation}\label{w1}
x=Re\int\Phi, \text{ where }
\end{equation}
\begin{equation}\label{w2}
\Phi=(\Phi_1,\Phi_2,\Phi_3)=((g^{-1}-g)\frac{dh}{2}, i(g^{-1}+g)\frac{dh}{2},dh).
\end{equation}
Conversely, let $M$ be a Riemann surface, $g\colon M \to \mathbb{C}\cup{\infty}$ a meromorphic function
and $dh$ a holomorphic one-form on $M$. Then, provided that $Re\int_{\alpha}\Phi=0$ for all closed curves $\alpha$ on $M$,  equation \eqref{w1} and \eqref{w2} define a conformal minimal
mapping of $M$ into $\rth$, which is regular provided the poles and zeros of
$g$ coincide with the zeros of $dh$. The holomorphic function $g$
and holomorphic one form $dh$ are the so-called {\it Weierstrass data}.
\end{proposition}

With Theorem~\ref{wt} in mind, it is easy to check that, for an immersed minimal surface $M$, admitting a nontrivial isometric deformation through minimal surfaces is equivalent to the condition that $Im\int_{\alpha}\Phi=0$ for all closed curves $\alpha$ on $M$. Take $p_1,...,p_n \in \C$ and, for any $k\in \N$, consider the following Weierstrass data 
$$g=\prod_{i=1}^n  (z-p_i)^{2k} \text{ and } dh=dz.$$
One can check that this data yield non-simply-connected, genus zero minimal surfaces admitting a nontrivial isometric deformation through minimal surfaces.

There are many results on the isometric indeformability of a constant mean curvature surface with topology (e.g. \cite{cmw1, ku2, me6, mr3, mt2, perez4}). In what follows, we give a criterion which guarantees the isometric indeformability of a constant mean curvature surface. In particular this result is a generalization of a rigidity theorem of Choi-Meeks-White for minimal surfaces, see Theorem 1.2 in~\cite{cmw1}.

\begin{theorem}\label{appl}
Let $x \colon M\to \rth$ be an isometric immersion of a smooth oriented surface with constant mean curvature $H$. Suppose that a  plane $\pi$ intersects $x(M)$ transversally in a closed unit speed curve $\gamma \colon [0,L] \to M$ then the component of the force of $\gamma$ orthogonal to $\pi$ is
$$\langle W([\gamma]), V\rangle =\int_0^L\frac{1}{\langle x_*(J\dot{\gamma}), V \rangle}\{\langle x_*(J\dot{\gamma}),V\rangle^2 +H\langle x\circ\gamma,\xi\rangle\} ds,$$
where $V$ is a unit vector normal to $\pi$.
If $$\int_0^L\frac{1}{\langle x_*(J\dot{\gamma}), V \rangle}\{\langle x_*(J\dot{\gamma}),V\rangle^2 +H\langle x\circ \gamma,\xi\rangle\} ds\neq 0 $$ then $x$ does not admit a nontrivial isometric deformation through surfaces of constant mean curvature $H$. In particular, if $x$ is minimal and the plane intersects it transversally in a closed curve then, to within congruences, $M$ admits only one other minimal isometric immersion in $\rth$, namely $-x$.
\end{theorem}

\begin{proof}
We may assume that $\pi$ is the $xy$-plane and $V=e_3$. Clearly, $$x\circ\gamma= \langle x\circ \gamma,x_*(\dot{\gamma})\rangle  x_*(\dot{\gamma})+\langle x\circ \gamma,x_*(J\dot{\gamma})\rangle x_*(J\dot{\gamma})+\langle x\circ \gamma,\xi\rangle \xi $$
and, if the origin is chosen in $\pi$,
$$0=\langle x\circ \gamma, e_3 \rangle= \langle x\circ \gamma,x_*(J\dot{\gamma})\rangle \langle x_*(J\dot{\gamma}),e_3\rangle+\langle x\circ \gamma,\xi \rangle \langle \xi,e_3\rangle.$$
Since the plane intersects $x(M)$ transversally $\langle x_*(J\dot{\gamma}),e_3\rangle \neq 0$ and consequently $$\langle x\circ \gamma,x_*(J\dot{\gamma})\rangle=-\frac{\langle x\circ \gamma,\xi\rangle \langle \xi, e_3 \rangle}{\langle x_*(J\dot{\gamma}),e_3\rangle}.$$
The decomposition of $x\circ \gamma$ above yields the following computation for $\langle W([\gamma]), V \rangle $ where $W([\gamma])$ is the force form defined in \S~\ref{force}.
\begin{align*}
\langle W([\gamma]), V \rangle &=\langle (\xi +Hx  \circ\gamma) \times x_*(\dot{\gamma}),e_3 \rangle\\
&=\langle x_*(J\dot{\gamma}) -H\langle x\circ \gamma, x_*(J\dot{\gamma}) \rangle \xi +H\langle x\circ \gamma, \xi \rangle x_*(J\dot{\gamma}) ,e_3 \rangle\\
&= \langle x_*(J\dot{\gamma}),e_3\rangle +H\{-\langle x \circ \gamma, x_*(J\dot{\gamma}) \rangle \langle \xi,e_3\rangle +\langle x\circ \gamma, \xi \rangle \langle x_*(J\dot{\gamma}) ,e_3 \rangle\}\\
&= \langle x_*(J\dot{\gamma}),e_3\rangle +H\langle x\circ \gamma,\xi\rangle \left \{ \frac{\langle \xi, e_3 \rangle^2}{\langle x_*(J\dot{\gamma}),e_3\rangle} +\langle x_*(J\dot{\gamma}) ,e_3 \rangle \right\}\\
&=\frac{1}{\langle x_*(J\dot{\gamma}), e_3\rangle}\{\langle x_*(J\dot{\gamma}), e_3 \rangle^2+H\langle x\circ\gamma , \xi\rangle \}
\end{align*}
The formula for $\langle W([\gamma]), V\rangle$ in Theorem~\ref{appl} now follows. Thus in the minimal case the force of $\gamma$ is clearly non-zero and the result follows from Theorem~\ref{main} and  Remark \ref{minimal}.

\end{proof}

\begin{corollary}
Let $x \colon M\to \rth$ be an isometric immersion of a smooth oriented surface with constant mean curvature $H\neq 0$. Suppose that $M$ has a  plane of symmetry $\pi$ which intersects $x(M)$ in a closed curve $\gamma \colon [0,L] \to M$. If $x\circ \gamma$ lies in an open disk of radius $\frac{1}{H}$ then there exists only finitely many isometric immersions of $M$ with constant mean curvature $H$.
\end{corollary}

\begin{proof}
We may assume that $\pi$ is the $xy$-plane and $V=e_3$. Since $\pi$ is a plane of symmetry, we can assume that $\langle x_*(J\dot{\gamma}), e_3 \rangle=1$. Therefore Theorem~\ref{appl} gives $$\langle W([\gamma]), e_3 \rangle=L+H\int_0^L \langle x\circ \gamma,\xi\rangle ds.$$

Suppose $x\circ \gamma$ lies in an open disk of radius $\frac{1}{H}$. Without loss of generality we can assume that the disk is centered at the origin. Then $|H\int_0^L \langle x\circ \gamma,\xi\rangle ds|<L$ and therefore $W([\gamma])\neq 0$. The result now follows from Theorem~\ref{main}.
\end{proof}

\begin{definition}
An embedding $x\colon M\to \rth$ has a plane of Alexandrov symmetry $\pi$ if $\pi$ is a plane of symmetry for $x(M)$ and $x(M)\backslash \{x(M)\cap \pi \}$ consists of two graphs over $\pi$.
\end{definition}

For a complete proper surface $x\colon M \to \rth$ with constant mean curvature $H\neq 0$ and a plane of Alexandrov symmetry Meeks-Tinaglia~\cite{mt4} showed that if $\sup_M|A|<\infty$ and $M$ has more then one end then $x$ has no other associates.
They show that $M$ has an end asymptotic to an unduloid (see also~\cite{kks1}). Since for unduloids none of the associates exists as we remarked in \S~\ref{force}, a compactness argument (which we also use at the end of the next proof) implies that  $x$ is the only isometric embedding of $M$ with constant mean curvature $H$ (see also~\cite{ku2,mt2}).

We prove the following theorem.

\begin{theorem}\label{alemb}
Let $x \colon M\to \rth$ be a complete proper isometric embedding of a smooth oriented surface with constant mean curvature $H\neq 0$, with a plane of Alexandrov symmetry. Then
\begin{itemize}
\item[(i)] there exists only finitely many isometric immersions of $M$ with constant mean curvature $H$, and
\item[(ii)] if $\sup_M|A|=\infty$ and $M$ has finite genus then $x$ and possibly $x_\pi$ are the only isometric embeddings of $M$ with constant mean curvature $H$. 
\end{itemize}
\end{theorem}

\begin{proof}
We begin by proving item {\it (i)}. Our goal is to show that there exists a plane which intersects $M$ transversally in a simple closed curve which is not homologically trivial. Once this is done, by a result of Korevaar and Kusner (Theorem 1.12 in~\cite{kk2}) the force of such a curve is always nonzero. Thus, by Theorem~\ref{main} there exist only finitely many isometric immersions of $M$ into $\rth$ with constant mean curvature $H$ and {\it (i)} will be proved.

Let $\Delta$ be the connected region of $\pi$ over which $x(M)\backslash \{x(M)\cap \pi \}$ is graphical. If the boundary of $\Delta$ consists of more then one connected component, select two components $\alpha$ and $\beta$. Let $\pi'$ be a plane perpendicular to $\pi$ which intersects both $\alpha$ and $\beta$ transversally. Choosing different boundary components of $\Delta$, if necessary, we can assume that the intersection of $\pi'$ with $\Delta$ contains a line segment connecting $\alpha$ and $\beta$. Since $x(M)$ is graphical over $\Delta$, $\pi'$ intersects $M$ transversally in a simple closed curve $\gamma$ intersecting $\alpha$ and $\beta$. Since $\alpha$ and $\beta$ are different components of $\partial \Delta$, $\gamma$  cannot be homologically trivial in $M$ and thus its force is nonzero by~\cite{kk2}.

If $\pi$ intersects $M$ in a single connected curve $\gamma$ and $\gamma$ is not closed, then $M$ would have exactly one end and genus zero. However, a result of Meeks says that a properly embedded surface with nonzero constant mean curvature and finite genus must have more than one end (see \cite{me17}). Thus $\gamma$ must be closed. 

If $\Delta$ is the unbounded region of $\pi$ outside $\gamma$ then $\gamma$ is homologically nontrivial. If $\Delta$ is the compact region of $\pi$ bounded by $\gamma$ then $x(M)\backslash \{x(M)\cap \pi \}$ consists of two compact graphs over $\Delta$ and $M$ must be compact. This implies that $x(M)$ is a round sphere in which case the isometric immersion is always unique. In either case, whether $\Delta$ is compact or not, there exist only finitely many isometric immersions of $M$ into $\rth$ with constant mean curvature $H$.

We now prove item {\it (ii)}. If $M$ has unbounded second fundamental form, then for any $n\in \mathbb{N}$ there exists $p_n\in M$ such that $|A|(p_n)>n$. Recall that for graphs with constant mean curvature and zero boundary value there exists a constant $C$ depending only on $H$ such that $\sup_M|x_3||A|<C$  (see for instance~\cite{ror1}) therefore in our case $|x_3(p_n)|<\frac{C}{|A|(p_n)}$ (here $x_3$ is the third component of the point $x$). After a sequence of translations which take $x(p_n)$ to the origin, we obtain a sequence of immersions $y_n\colon M \to \rth$, $y_n:=x-x(p_n)$, with constant mean curvature and a plane of Alexandrov symmetry such that $y_n(p_n)=0$. The distance from the origin to the plane of symmetry of $y_n$ is bounded by  $\frac{C}{|A|(p_n)}$.

Consider the sequence of non-negative
functions $F_n \colon M \to \mathbb{R}$,
$$F_n(p)=(|y_n(p)|-1)^2|A|^2(p)$$ over the connected component $M_n$ of $\{p\in M \colon |y_n(p)|\leq 1 \}$ containing $p_n$. 
The function $F_n$ is zero on the boundary of
$M_n$ and therefore it attains its maximum on $M_n$ at a point in its interior.
Let $q_n$ be such a point, i.e.
$$F_n(q_n)=(|y_n(q_n)|-1)^2|A|^2(q_n)=\max_{M_n}F_n(p)\geq F_n(p_n)=|A|^2(p_n).$$
Fix $\sigma_n>0$ such that $2\sigma_n<1-|y_n(q_n)|$ and
$$4\sigma_n^2|A|^2(q_n)=4|A|^2(p_n)=C^2_n.$$ Notice that $\sigma_n$ is less than half the distance from $y_n(q_n)$ to the boundary of the ball of radius one centered at the origin. Let $M_{\sigma_n}$ be the connected component of $\{p\in M_n \colon |y_n(p) - y_n(q_n)|\leq \sigma_n \}$ containing $q_n$. Since $F_n$ achieves its maximum on $M_n$ at
$q_n$,
\begin{align*}\sup_{M_{\sigma_n}}\sigma_n^2|A|^2 &\leq
\sup_{M_{\sigma_n}}\sigma_n^2 \frac{F_n(p)}{(|y_n(p)|-1)^2}\\
& \leq\frac{4\sigma_n^2}{(|y_n(q_n)|-1)^2}\sup_{M_{\sigma_n}}F_n(p)\\
& = \frac{4\sigma_n^2}{(|y_n(q_n)|-1)^2}F_n(q_n) =4\sigma_n^2|A|^2(q_n).
\end{align*}
For any $n$, we apply a translation which takes the plane of Alexandrov symmetry to the $xy$-plane and $y_n(q_n)$ on the $z$-axis, and let $z_n\colon M\to \rth$, $z_n:=y_n+v_n$ for a certain $v_n\in \rth$, denote this new sequence of immersions. We have obtained the following
$$\sup_{M_{\sigma_n}}|A|^2\leq 4|A|^2(q_n),\, 4\sigma_n^2|A|^2(q_n)=C_n^2, \text{ and } |(z_n)_3(q_n)|<\frac{C}{|A|(q_n)}.$$ 

Consider a new sequence of immersions $w_n\colon M \to \rth$ obtained by rescaling $z_n$ by a factor of $|A|(q_n)$, $w_n:=|A|(q_n)z_n$. Note that $|(w_n)_3(q_n)|<C$ and that $|H_n|\leq \frac{|H|}{|A|(q_n)}$, where $H_n$ denotes the mean curvature of the immersion $w_n$. Since we are assuming that the genus of $M$ is finite, a standard compactness argument implies that this sequence converges in the $C^2$ convergence to a non-flat minimal embedding of a genus zero surface, $x_\infty \colon M_\infty \to \rth$, with bounded second fundamental form and hence properly embedded (see~\cite{cm35} and also~\cite{mr13}). Since a proper embedding cannot be contained in a half-space (see~\cite{hm10}), the $xy$-plane must be a plane of symmetry. This implies that $x_\infty\colon M_\infty \to \rth$ must be a catenoid (see~\cite{mpr6} and \cite{cm25,cm21,cm22,cm24,cm23,col1,lor1,mr8} and others). The following compactness argument (see also~\cite{mt2}) then implies that in either case, $x$ and possibly $x_\pi$ are the only isometric embeddings of $M$ into $\rth$ with constant mean curvature~$H$. 

Let $\gamma_\infty \colon [0,L] \to M_\infty$ denote the shortest closed geodesic on $M_\infty$  and let $\gamma_n \colon [0,L_n] \to M_{\sigma_n}$, be the sequence of cycles in $M_{\sigma_n}$ such that $z_n(\gamma_n)$ converges to $x_\infty(\gamma_\infty)$. Suppose that there exists $\theta$ different from zero or $\pi$ such that the associate $x_\theta \colon M \to \rth$ exists then, because of the convergence, 
$(z_n)_\theta(\gamma_n)$ converge to $(x_\infty)_\theta(\gamma_\infty)$. However, while $(z_{n})_\theta(\gamma_n)$ must be a closed curve, an easy computation shows that when $\theta$ is different from zero or $\pi$, $(x_{\infty})_\theta(\gamma_\infty)$ is not a closed curve. This contradicts the convergence and proves that $x$ and possibly $x_\pi$ are the only isometric embeddings of $M$ into $\rth$ with constant mean curvature~$H$.


\end{proof}

\bibliographystyle{plain}
\bibliography{bill}

\end{document}

%% file: orientation.pstex_t
\begin{picture}(0,0)%
\includegraphics{orientation.pstex}%
\end{picture}%
\setlength{\unitlength}{3947sp}%
\begingroup\makeatletter\ifx\SetFigFont\undefined%
\gdef\SetFigFont#1#2#3#4#5{%
  \reset@font\fontsize{#1}{#2pt}%
  \fontfamily{#3}\fontseries{#4}\fontshape{#5}%
  \selectfont}%
\fi\endgroup%
\begin{picture}(3117,2305)(4864,-4955)
\put(5026,-3961){\makebox(0,0)[lb]{\smash{{\SetFigFont{12}{14.4}{\familydefault}{\mddefault}{\updefault}{\color[rgb]{0,0,0}D}%
}}}}
\put(5401,-3061){\makebox(0,0)[lb]{\smash{{\SetFigFont{12}{14.4}{\familydefault}{\mddefault}{\updefault}{\color[rgb]{0,0,0}F}%
}}}}
\put(6151,-3586){\makebox(0,0)[lb]{\smash{{\SetFigFont{12}{14.4}{\familydefault}{\mddefault}{\updefault}{\color[rgb]{0,0,0}$\xi$}%
}}}}
\put(6826,-4036){\makebox(0,0)[lb]{\smash{{\SetFigFont{12}{14.4}{\familydefault}{\mddefault}{\updefault}{\color[rgb]{0,0,0}$\eta$}%
}}}}
\put(7501,-4036){\makebox(0,0)[lb]{\smash{{\SetFigFont{12}{14.4}{\familydefault}{\mddefault}{\updefault}{\color[rgb]{0,0,0}K}%
}}}}
\put(7726,-3511){\makebox(0,0)[lb]{\smash{{\SetFigFont{12}{14.4}{\familydefault}{\mddefault}{\updefault}{\color[rgb]{0,0,0}$\gamma$}%
}}}}
\put(7651,-4336){\makebox(0,0)[lb]{\smash{{\SetFigFont{12}{14.4}{\familydefault}{\mddefault}{\updefault}{\color[rgb]{0,0,0}$\nu_k$}%
}}}}
\put(6901,-4411){\makebox(0,0)[lb]{\smash{{\SetFigFont{12}{14.4}{\familydefault}{\mddefault}{\updefault}{\color[rgb]{0,0,0}$\dot{\gamma}$}%
}}}}
\end{picture}%